\documentclass[11pt]{article}
\usepackage{amsmath}
\usepackage{mathrsfs}
\usepackage{amsfonts}

\usepackage{amsfonts, amsmath, amssymb}
\usepackage{amssymb,amsfonts,amsmath,color,
latexsym, epsfig,cite, psfrag,eepic,colordvi}
\usepackage{amscd,graphics}
\usepackage{lineno}
\usepackage{comment}

\newcommand{\qed}{\hfill $\Box $}
\newcommand{\pf}{\noindent {\bf Proof.} }

\newtheorem{theorem}{Theorem}[section]
\newtheorem{lemma}[theorem]{Lemma}

\begin{document}

\title{A Degree Condition for a Graph to have  $(a,b)$-Parity Factors \thanks{This work is supported   by  the National Natural Science Foundation
of China (Grant No. 11471257)}}

\author{Haodong Liu and Hongliang Lu\thanks{Corresponding email: luhongliang215@sina.com (H. Lu)}\\ {\small \textsuperscript{a}Department of Mathematics}
\\ {\small Xi'an Jiaotong University, Xi'an 710049, PR China}
}

\date{}

\maketitle

\begin{abstract}
Let $a,b,n$ be three positive integers such that $a\equiv b\pmod 2$ and $n\geq b(a+b)(a+b+2)/(2a)$.
Let $G$ be a graph of order $n$ with minimum degree at least $a+b/a-1$. We show that $G$ has an $(a,b)$-parity factor, if $max\{d_G(u),d_G(v)\}\geq \frac{an}{a+b}$ for any two  nonadjacent vertices  $u,v$ of $G$. It is an extension of Nishimura's results for the existence of $k$-factors (\emph{J. Graph Theory}, \textbf{16} (1992), 141--151) and generalizes Li and Cai's result in some senses (\emph{J. Graph Theory}, \textbf{27} (1998),  1--6). These  conditions are tight.
\end{abstract}

\noindent\textbf{Keywords: degree condition, parity factor}

\section{Introduction}

In this paper we consider only simple graphs. Let $G$ be a graph with vertex set $V(G)$ and
edge set $E(G)$.
Given   $v\in V(G)$, let $N_G(v)$ denote the set of vertices adjacent to with $v$ in $G$ and $d_G(v)=|N_G(v)|$. The minimum vertex degree in graph $G$ is denoted by $\delta(G)$.
We write
$N_G[v]=N_G(v)\cup \{v\}$. Given $D\subseteq V(G)$, let $N_D(V)=N_G[v]\cap D$. For $X\subseteq V(G)$, the subgraph of $G$ whose vertex set is $X$ and whose edge set consists of the edges of $G$ joining vertices of $X$ is called the subgraph of $G$ induced by $X$ and is denoted by $G[X]$.

Let $g,f$ be two non-negative integer-valued function such that $g(v)\leq f(v)$ and $g(v)\equiv f(v)\pmod 2$ for all $v\in V(G)$.
A spanning subgraph $F$ of $G$ is called \emph{$(g,f)$-parity factor} if $d_F(v)\equiv f(v)\pmod 2$ and $g(v)\leq d_F(v)\leq f(v)$ for all $v\in V(G)$.
A $(g,f)$-parity factor is called \emph{$f$-factor} if $f(v)=g(v)$ for all $v\in V(G)$. If $f(v)=k$ for all $v\in V(G)$, then an $f$-factor is called a \emph{$k$-factor}. Let $a,b$ be two integers such that $a\leq b$ and $a\equiv b\pmod 2$. If $f(v)=b$ and $g(v)=a$ for all $v\in V(G)$, then a $(g,f)$-parity factor is called an \emph{$(a,b)$-parity factor}.

Lov\'asz \cite{Lov72} gave a characterization of graphs having $(g,f)$-parity factors.   Amahashi \cite{Amha85} found a Tutte's type characterization for $(1,k)$-odd factors, which was  generalized  to $(1,f)$-odd factors by  Cui and Kano \cite{CK88}.
\begin{theorem}[Lov\'asz, \cite{Lov72}]\label{lov72}
A graph $G$ has a $(g,f)$-parity factor if and only if for any two disjoint subsets $S,T$ of $V(G)$,
\[
\eta(S,T)=f(S)-g(T)+\sum_{x\in T}d_{G-S}(x)-q(S,T)\geq 0,
\]
where $q(S,T)$ denotes the number of  components $C$ of $G-S-T$, called $g$-odd components, such that $g(V(C))+e_G(V(C),T)\equiv 1 \pmod 2$.
\end{theorem}

\begin{theorem}[Amahashi, \cite{Amha85}]\label{Am85}
Let $k\geq 1$ be an odd integer.
A graph $G$ contains an $(1,k)$-parity factor if and only if for any subset $S\subseteq V(G)$,
\[
c_0(G-S)\leq k|S|,
\]
 where
$c_0(G-S)$ denotes the number of odd components of $G-S$.
\end{theorem}

Nishimura{\cite{Nishi}} gave a degree conditions for a graph to have a $k$-factor.

\begin{theorem}\label{Ni92}
Let $k$ be an integer such that $k\geq 3$, and let $G$ be a connected
graph of order $n$ with $n\geq 4k-3$, $kn$ even, and minimum degree at
least $k$. If $G$ satisfies
\begin{equation*}
\max\{d_G(u),d_G(v)\}\geq n/2
\end{equation*}
for each pair of nonadjacent vertices $u,v$ in $G$, then $G$ has a  $k$-factor.
\end{theorem}

Li and Cai  {\cite{Cai}} give a degree conditions for a graph to have an $(a,b)$-factor, which extended Nishimura's result.

\begin{theorem}
Let $G$ be a graph of order $n$, and let $a$ and $b$ be integers such that $1\leq a<b$.
Then $G$ has an [a,b]-factor if $\delta(G)\geq a, n\geq 2a+b+\frac{a^2-a}{b}$ and
\begin{align}\label{main-degree-condition}
  \max\{d_G(u),d_G(v)\}\geq \frac{an}{a+b}
\end{align}
for any two nonadjacent vertices $u$ and $v$ in $G$.
\end{theorem}

In this paper we give a sufficient condition for a graph to have an $(a,b)$-parity factor in term of the minimum degree of graph $G$.  Our main result generalizes Nishimura's result and improves Li and Cai's  result in some sense.

\begin{theorem}\label{Main-Theorem}
Let $a,b,n$ be three integers such that $a\equiv b\pmod 2$, $na$ is even and $n\geq b(a+b)(a+b+2)/(2a)$.
Let $G$ be a graph of order $n$. If  $\delta (G)\geq a+\frac{b-a}{a}$ and
\begin{equation}\label{degree_condition}
max\{d_G(u),d_G(v)\}\geq \frac{an}{a+b}
\end{equation}
for any two nonadjacent vertices, then $G$ has an $(a,b)$-parity factor.
\end{theorem}

\section{Proof of Theorem \ref{Main-Theorem}}

Firstly, we show that Theorem \ref{Main-Theorem} holds for $a=1$.
\begin{lemma}\label{a=1}
Let $k,n$ be two positive integers such that $k$ is odd, $n$ is even and $n\geq k+1$.
Let $G$ be a connected graph with  order $n$. If $G$ satisfies
\begin{equation}\label{degree_condition2}
max\{d_G(u),d_G(v)\}\geq \frac{n}{1+k}
\end{equation}
for each pair of nonadjacent vertices, then $G$ has a $(1,k)$-odd factor.
 \end{lemma}

 \pf Suppose that $G$ contains no $(1,k)$-parity factors. By Theorem \ref{Am85}, there exists a subset $S\subset V(G)$ such that
 \begin{align*}
 c_o(G-S)> k|S|.
 \end{align*}
 Let $C_1,\ldots, C_q$ be these odd components of $G-S$ such that $|C_1|\leq \cdots\leq |C_q|$. Note that $n$ is even and $G$ is connected. By parity, one can see that
\begin{align}\label{q>=ks+2}
 q=c_o(G-S)\geq k|S|+2\ \mbox{and}\ S\neq \emptyset.
\end{align}
 Let $u\in V(C_1)$ and $v\in V(C_2)$.
 By (\ref{q>=ks+2}), we infer that
 \begin{align*}
 \max\{d_G(u),d_G(v)\}&\leq \{|C_1|-1+|S|,|C_2|-1+|S|\}\\
 &\leq |C_2|-1+|S|\\
 &\leq \frac{n-|S|-1}{q-1}-1+|S|\\
 &\leq \frac{n}{k|S|+1}+|S|-1-\frac{1}{k}\quad \mbox{(since $|S|\geq 1$ and $n\geq (k+1)|S|+2$)}\\
 &< \frac{n}{k+1},
 \end{align*}
 contradicts with (\ref{degree_condition2}) since $uv\notin E(G)$. \qed

\noindent\textbf{Proof of Theorem \ref{Main-Theorem}.}
By Lemma \ref{a=1}, we may assume that $a\geq 2$. By Theorem \ref{Ni92}, we may assume that $a\leq b-2$.
Suppose that $G$ contains no $(a,b)$-parity factors. By Theorem \ref{lov72}, there exists two disjoint vertex sets $S$ and $T$ such that
\begin{align}\label{eq_lov1}
\eta(S,T)=b|S|-a|T|+\sum_{x\in T}d_{G-S}(x)-q(S,T)\leq -2,
\end{align}
where $q(S,T)$ denotes the number of  components $C$ of $G-S-T$, called $a$-odd components, such that $g(V(C))+e_G(V(C),T)\equiv 1 \pmod 2$.
We write $s=|S|$, $t=|T|$ and $w=q(S,T)$. From (\ref{eq_lov1}), one can see that
\begin{align}\label{eq_lov2}
\eta(S,T)=bs-at+\sum_{x\in T}d_{G-S}(x)-w\leq -2.
\end{align}

If $S\cup T=\emptyset$, we have $w\geq 2$ by (\ref{lov72}), which implies that $G$ consists of
at least two components. However, this contradicts the connectedness of $G$.
So we may assume that
\begin{align}\label{ST-nonempty}
S\cup T\neq \emptyset.
\end{align}
If $w\geq 1$, let
$C_1,C_2,\dots C_w$ denote these $a$-odd components of $G-S-T$, and $m_i=|V(C_i)|$ for $1\leq i\leq w$. Put $U=\bigcup_{1\leq i\leq w}V(C_i)$.

We pick $S$ and $T$ such that $U$ is minimal and $V(G)-S-T-U$ is maximal.

\medskip
\textbf{ Claim 1.~}   $d_{G-S}(u) \geq a+1$  and $ e_{G}(u,T) \leq b-1$ for every vertex  $u\in V(U)$.
\medskip

Firstly, suppose that there exists $u\in U$ such that
\[
d_{G-S}(u) \leq a.
\]
Let $T'=T\cup \{u\}$.  One can see that
\begin{align*}
\eta(S,T')&=bs-a|T'|+\sum_{x\in T'}d_{G-S}(x)-q(S,T')\\
&=bs-at-a+\sum_{x\in T}d_{G-S}(x)+d_{G-S}(u)-q(S,T')\\
&\leq bs-at+\sum_{x\in T}d_{G-S}(x)-(q(S,T)-1)\\
&\leq -1,
\end{align*}
which implies by parity
\begin{align*}
\eta(S,T')=bs-a|T'|+\sum_{x\in T'}d_{G-S}(x)-q(S,T')\leq -2,
\end{align*}
contradicting the minimality of $U$.
Secondly, suppose that there exists $u\in U$ such that
\[
e_G(u,T) \geq b.
\]
Let $S'=S\cup \{u\}$.  One can see that
\begin{align*}
\eta(S',T)&=b|S'|-at+\sum_{x\in T}d_{G-S'}(x)-q(S',T)\\
&=bs+b-at+\sum_{x\in T}d_{G-S}(x)-e_G(u,T)-q(S',T)\\
&\leq bs-at+\sum_{x\in T}d_{G-S}(x)-(q(S,T)-1)\\
&\leq -1,
\end{align*}
which implies by parity
\begin{align*}
\eta(S',T)=b|S'|-at+\sum_{x\in T}d_{G-S'}(x)-q(S',T)\leq -2,
\end{align*}
contradicting the minimality of $U$ again.  This completes Claim 1. \qed

\medskip
\textbf{ Claim 2.~} Let $C_{i_1},\ldots,C_{i_{\tau}}$ be any $\tau$ components of $G[U]$ and let $U'=\bigcup_{j=1}^{\tau}V(C_{i_j})$.  $d_{G[T\cup U']}(u)\leq  a-1+\tau$ for every vertex  $u\in T$.
\medskip

 Suppose that there exists $u\in T$ such that $d_{G[T\cup U']}(u)\geq  a+\tau$. Let $T'=T-u$.
One may see that
\begin{align*}
\eta(S,T')&=bs-a|T'|+\sum_{x\in T'}d_{G-S}(x)-q(S,T')\\
&= bs-at+a+\sum_{x\in T}d_{G-S}(x)-d_{G-S}(u)-q(S,T')\\
&\leq bs-at+a+\sum_{x\in T}d_{G-S}(x)-(a+\tau)-(q(S,T)-\tau)\\
&=bs-at+\sum_{x\in T}d_{G-S}(x)-q(S,T)\leq -2,
\end{align*}
contradicting to the maximality of $V(G)-S-T-U$. This completes Claim 2. \qed

From the definition of $U$, we have
\begin{equation}\label{maxorder}
  |U| \geq m_1+m_2(w-1).
\end{equation}
By Claim 1, one can see that for every $u\in C_j\ (1\leq j\leq w)$,
\begin{equation}\label{eq-m_j-s-r}
  d_G(u)\leq (m_{j}-1)+s+r
\end{equation}
where $r=\min\{b,t\}$.
Let $u_1\in V(C_1)$ and $u_2\in V(C_2)$. It follow from   (\ref{eq-m_j-s-r}) that
\begin{equation}\label{eqm1m2}
    \max\{d_G(u_1),d_G(u_2)\}\leq (m_{2}-1)+s+r
\end{equation}

\medskip
\textbf{Claim 3.} $S\neq\emptyset$.
\medskip

Suppose that $S=\emptyset$.
By (\ref{ST-nonempty}), one may see that  $t\geq 1$.
Note that $\delta(G)\geq a+\frac{b-a}{a}$. So we have
 $d_G(x_1)\geq a+\frac{b-a}{a}$. From Theorem (\ref{lov72}), one can see that
\begin{align}\label{w-low-bound}
  w=q(S,T)\geq \sum\limits_{v\in T}d_{G}(v)-at+2 \geq  \frac{b-a}{a}t+2.
\end{align}
If $w>\frac{b}{a}+2$, since $b\geq a+2$, then it follows that
\begin{align}\label{w>=(b+2)/a+2}
w\geq \frac{b+2}{a}+2.
\end{align}
Combining (\ref{eqm1m2}) and (\ref{w>=(b+2)/a+2}), one can see that
\[
 \max\{d_G(u_1),d_G(u_2)\}\leq m_2-1+s+r\leq \frac{n-t-1}{w-1}+b<\frac{an}{a+b+2}+b<\frac{an}{a+b},
\]
contradicting to (\ref{main-degree-condition}).
So we may assume that $w\leq \frac{b}{a}+2$. From (\ref{w-low-bound}), we infer that
\[
\frac{b-a}{a}t+2\leq \frac{b}{a}+2,
\]
i.e.,
\begin{align}\label{t-bound-cl3}
t\leq \frac{b}{b-a}.
\end{align}
From (\ref{w-low-bound}), we have
\begin{align*}\label{m1_cl1}
m_1\leq \frac{a(n-t)}{a+b}.
\end{align*}
Consider $H=G[V(C_1)\cup T]$. By Claim 2, for every $y\in T$, one can see that
\[
d_G(y)\leq a-1+w\leq \frac{b}{a}+1+a<\frac{an}{a+b}.
\]
By Claims 1 and 2, $d_{G-S}(u)=d_H(u)\geq a+1$ for every $u\in V(C_1)$ and $d_H(v)\leq a$ for every $v\in T$. Thus there exists two non-adjacent vertices $u\in V(C_1)$ and $v\in T$. If $m_1\leq \frac{an}{a+b}-\frac{b}{b-a}$, then one can see that
\[
\max\{d_G(u),d_G(v)\}\leq m_1-1+t<\frac{an}{a+b},
\]
a contradiction.
Thus we may assume that $ m_1\geq \frac{an}{a+b}-\frac{b}{b-a}$.
We claim that there exists $u\in V(C_1)$ such that $e_G(u,T)=0$, otherwise, by Claim 2 and (\ref{t-bound-cl3}), we have
\[
\frac{ab}{b-a}\geq at\geq \sum_{x\in T}d_{G[V(C_1)\cup T]}(x)\geq m_1\geq \frac{an}{a+b}-\frac{b}{b-a},
\]
i.e.,
\[
n\leq \frac{(a+1)b(a+b)}{a(b-a)},
\]
a contradiction. It follows
\[
\max\{d_G(u),d_G(v)\}\leq m)\leq m_1-1+s=m_1-1<\frac{an}{b+a},
\]
a contradiction.  This completes Claim 3. \qed

\medskip
\textbf{Claim 4.} $T\neq\emptyset$.
\medskip

Suppose that $T=\emptyset$. By Theorem \ref{lov72},  then we have
\[
w\geq bs+2.
\]
By Claim 1, we have $|V(C_i)|\geq a+1$ for $1\leq i\leq w$. Thus we infer that
\[
n\geq (a+1)w+s\geq (a+1)(bs+2)+s>(a+1)bs,
\]
which implies that
\begin{equation*}
  s\leq \frac{n}{(a+1)b}.
\end{equation*}
Hence we have
\[
m_2-1+s\leq \frac{n-s}{w-1}+s\leq \frac{n}{bs+1}+s<\frac{an}{a+b},
\]
a contradiction. This completes Claim 4. \qed

Put $h_1:=\min\{ d_{G- S}(v)\ |\ v\in T\}$,
and let $x_1\in T$ be a vertex satisfying $d_{G- S}(x_1)=h_1$. We write $p=|N_T[x_1]|$. Further, if $T- N_T[x_1]\neq\emptyset$,
let $  h_2:=\min\{ d_{G- S}(v)\ |\ v\in T-N_T[x_1]\}$ and let $x_2\in T-N_T[x_1]$ such that $d_{G- S}(x_2)=h_2$.
By the definition of $x_i$, we have
\begin{align}\label{dx1x2}
 \max\{d_G(x_1),d_G(x_2)\}&\leq  \max\{h_1+s,h_2+s\}\leq h_2+s.
\end{align}

Now we discuss four cases.

\medskip
\textbf{Case 1.} $h_1\geq a$.
\medskip

By Theorem \ref{lov72}, one can see that
\begin{align*}
  w &\geq bs-at+\sum_{v\in T}d_{G- S}(v)+2\\
  & \geq bs+(h_1-a)t+2\\
  &\geq bs+2,
\end{align*}
i.e.,
\begin{align}\label{lov-case1-1}
w \geq bs+2.
\end{align}
Note that $n\geq w+s+t$. From (\ref{lov-case1-1}), we infer that
\begin{align*}
s<\frac{n}{b+1}.
\end{align*}
Hence we have
\begin{align*}
m_2-1+s+r&\leq \frac{n-t}{bs+1}+s+b\\
&\leq \frac{n}{bs+1}+s+b-1\\
&<\frac{an}{a+b},
\end{align*}
a contradiction.

So we  may assume that $h_1<a$.

\medskip
\textbf{Case 2.} $T=N_T[x_1]$.
\medskip

We write  $t=|N_T[x_1]|$. Since $h_1<a$, we have  $t\leq a$. By Claim 1, one can see that for every $u\in V(C_1)$, $d_{G-S}(u)\geq a+1>h_1$. Thus we infer that $V(C_1)-N_G(x_1)\neq \emptyset$, i.e.,  there exists a vertex $v\in V(C_1)$ such that $x_1v\notin E(G)$.
By Theorem \ref{lov72},
\begin{align*}
  w&\geq bs+(h_1-a)t+2\\
  &\geq bs+(h_1-a)(h_1+1)+2\\
  &\geq b(a+\frac{b}{a}-1-h_1)+(h_1-a)(h_1+1)+2\quad \mbox{(since $s+h_1\geq\delta(G)\geq a+\frac{b}{a}-1$)}\\
  &=h_1^2-(a+b-1)h_1+ab+\frac{b^2}{a}-a-b+2\quad \mbox{(since $1\leq h_1\leq a$ and $a< b$)}\\
  &\geq \frac{b^2}{a}-b+2>0,
\end{align*}
i.e.,
\begin{align}\label{Case21-w}
  w &\geq bs+(h_1-a)(h_1+1)+2>2.
\end{align}
One can see that
\begin{align*}
\max\{d_G(v),d_G(x_1)\}&\leq m_1+s+t-1\\
&\leq \frac{n-s-t}{w}+s+t-1\\
&\leq \frac{n-s-t}{bs+(h_1-a)(h_1+1)+2}+s+t-1\\
& \leq \frac{n-s-h_1-1}{bs+(h_1-a)(h_1+1)+2}+s+h_1\\
&=\frac{n-h_1-1+\frac{1}{b}(h_1-a)(h_1+1)}{bs+(h_1-a)(h_1+1)+2}-\frac{1}{b}+s+h_1,
\end{align*}
i.e.,
\begin{align*}
\max\{d_G(v),d_G(x_1)\}&\leq \frac{n-h_1-1+\frac{1}{b}(h_1-a)(h_1+1)}{bs+(h_1-a)(h_1+1)+2}-\frac{1}{b}+s+h_1.
\end{align*}
We write
\begin{align}
f(s)= \frac{n-h_1-1+\frac{1}{b}(h_1-a)(h_1+1)}{bs+(h_1-a)(h_1+1)+2}-\frac{1}{b}+s+h_1.
\end{align}
So we have
\begin{align}\label{f'(s)}
f'(s)=-\frac{b(n-h_1-1)+(h_1-a)(h_1+1)}{(bs+(h_1-a)(h_1+1)+2)^2}+1.
\end{align}

Now we discuss two subcases.

\medskip
\textbf{Case 2.1.~} $s\leq \frac{an}{a+b}-h_1-1$.
\medskip

By (\ref{Case21-w}) and (\ref{f'(s)}),  we infer that
\begin{align}\label{cov-f(s)}
f(s)\leq \max \{f(a+\frac{b}{a}-h_1-1), f(\frac{an}{a+b}-h_1-1)\}.
\end{align}
Hence one can see that
\begin{align*}
max\{d_G(x_1),d_G(x_2)\}&\leq m_1-1+s+t\\
&\leq \max\{f(a+\frac{b}{a}-h_1-1), f(\frac{an}{a+b}-h_1-1)\}\\
&< \frac{an}{a+b},
\end{align*}
contradicting to the degree condition.

\medskip
\textbf{Case 2.2.~} $s>\frac{an}{a+b}-h_1-1$.
\medskip

One can see that
\begin{align*}
n&\geq s+t+w\\
&\geq s+t+bs+(h_1-a)t+2\\
&\geq (b+1)\frac{an}{a+b}-(b+1)h_1-(b+1)+(h_1-a)(h_1+1)+2\\
&=(b+1)\frac{an}{a+b}+h_1^2-(a+b)h_1-a-b+1\quad \mbox{(since $0\leq h_1\leq a$)}\\
&\geq \frac{abn}{a+b}+\frac{an}{a+b}-ab-a-b+1>n \quad \mbox{(since $a\geq 2$)},
\end{align*}
a contradiction.

\medskip
\textbf{Case 3.~} $h_2\geq a$.
\medskip

By Lov\'{a}sz Theorem \ref{lov72},
\begin{align*}
  w&\geq   bs+\sum\limits_{v\in T}d_{G-S}(v)-at+2 \\
  &\geq   bs+(h_1-a)p+(h_2-p)(t-p)+2 \\
  &\geq  bs+(h_1-a)p+2
\end{align*}

Now we discuss two subcases.

\medskip
\textbf{Subcase 3.1.~} $h_2\leq \frac{1}{4}(a^2+6a+5)$.
\medskip

By Lov\'asz Theorem \ref{lov72}, we find
$$s\geq \frac{an}{a+b}-h_2\geq\frac{an}{a+b}-\frac{1}{4}(a^2+6a+5).$$
Hence one can see that
\begin{align*}
n&\geq w+s+t\\
&\geq  (b+1)s+(h_1-a)p+2+t\\
&\geq (b+1)s+(h_1-a+1)p+2\\
&\geq (b+1)(\frac{an}{a+b}-\frac{1}{4}(a^2+6a+5)+1)+(h_1-a+1)(h_1+1)+2\\
&>n,
\end{align*}
a contradiction.

\medskip
\textbf{Subcase 3.2.~} $h_2\geq \frac{1}{4}(a^2+6a+5)$.
\medskip

By Lov\'{a}sz Theorem \ref{lov72},
\begin{align*}
  w&\geq   bs+\sum\limits_{v\in T}d_{G-S}(v)-at+2 \\
  &\geq   bs+(h_1-a)p+(h_2-a)(t-p)+2 \\
  &\geq  bs+(h_1-a)(h_1+1)+h_2-a+2\\
  &\geq bs-\frac{1}{4}(a+1)^2+h_2-a+2\\
  &\geq bs+3.
\end{align*}
We find
\begin{align*}
n\geq s+t+w\geq (b+1)s+2,
\end{align*}
which implies that
\begin{align*}
s\leq \frac{n-2}{b+1}.
\end{align*}
Thus by Claim 1, we infer that
\begin{align*}
m_2-1+s+b&\leq \frac{n-s-t}{bs+2}+s+b-1\\
&\leq \frac{n-s}{bs+2}+s+b-1\\
&<\frac{n-2}{bs+1}+s+b\\
&<\frac{an}{a+b},
\end{align*}
a contradiction.

\medskip
\textbf{Case 4.~} $0\leq  h_1\leq h_2\leq a-1$.
\medskip

By (\ref{main-degree-condition}), we infer that
\begin{align}\label{case4s-bound}
s\geq \frac{an}{a+b}-h_2.
\end{align}
By Lovasz Theorem \ref{lov72},
\begin{align*}
  w &\geq   bs+\sum\limits_{v\in T}d_{G-S}(v)-at+2 \\
  &\geq  bs+(h_1-a)p+(h_2-a)(t-p)+2,
\end{align*}
where $p=|N_T[x_1]|$ and naturally there is $p\leq a$,
\begin{align}\label{case4w-bound}
  w &\geq  bs+(h_1-a)p+(h_2-a)(t-p)+2.
\end{align}
Thus we get
\begin{align*}
  n&\geq s+t+w\\
  &\geq (b+1)s+(h_1-a)p+(h_2-a)(t-p)+2+t\\
  &= (b+1)s+(h_1-h_2)p+(h_2+1-a)t+2\\
  &\geq (b+1)(\frac{an}{a+b}-h_2)+(h_1-h_2)(h_1+1)+(h_2-a+1)(\frac{bn}{a+b}+h_2)+2\\
  &\geq (b+1)\frac{an}{a+b}+(-a+1)(\frac{bn}{a+b}+h_2)+2+h_2(\frac{bn}{a+b}+h_2-b-1)+(h_1-h_2)(h_1+1)\\
  &\geq n+2,
\end{align*}
a contradiction.
This completes the proof. \qed

\noindent\textbf{Remark :}
These minimum degree conditions are sharp. Let $a,b,m$ be three integers, such that $m$ is sufficiently large. Consider graph $K_{ma,mb+1}$. Denote $K_{ma,mb+1}$ by $G$.   Li  and Cai \cite{Cai} show that $K_{ma,mb+1}$ contains no $[a,b]$-factors. One can see that
\[
\frac{an}{a+b}>\delta(G)\geq ma> \frac{a|V(G)|}{a+b}-1.
\]
Other hand. Let $C_1,\ldots, C_q$ be $q$ copies of $K_m$, where $q=a+\lceil\frac{b-a}{a}\rceil-1$. Let $G'$ be a graph obtained form $C_1,\ldots,C_q$ by adding a new vertex $v$ connecting one of vertices of each copy. Clearly, $G'$ is connected and $\delta(G')= a+\lceil\frac{b-a}{a}\rceil-1$.  By taking $S=\emptyset$ and $T=\{v\}$, we infer that $G'$ contains no $(a,b)$-parity factor by   Lovasz's Theorem \ref{lov72}.

\end{document}